\documentclass [11pt] {article}

\usepackage{amsfonts}
\usepackage{amsmath}
\usepackage{amssymb}
\usepackage{color}
\usepackage{graphicx}
\usepackage{times}
\usepackage[round]{natbib}
\usepackage[T1]{fontenc}
\usepackage{url}

\newcommand {\sirev} {\textit {SIAM Review\/}}
\newcommand {\sandr} {\textit {Survey and Review\/}}
\newcommand {\pandt} {\textit {Problems and Techniques\/}}
\newcommand {\expo} {\textit {Expository Research Papers\/}}
\newcommand {\sigest} {\textit {SIGEST\/}}
\newcommand {\educ} {\textit {Education\/}}

\newcommand {\pands} {\textit {Problems and Solutions\/}}

\begin{document}

\title{A Brief Review of \textit {SIAM Review}}

\author{Joseph F. Grcar\,\footnote{6059 Castlebrook Drive, Castro Valley, CA 94552 USA.}\, \footnote {\texttt {jfgrcar@comcast.net}, \texttt {jfgrcar@gmail.com.}}}

\date{}

\maketitle

\begin{abstract}
\noindent
\sirev\ is examined for referee delay, citations, and paper length after the reorganization of the journal in 1999. A single, very-highly cited article was responsible for all the increase to the impact factor during the past decade; the reorganization did not improve the journal overall. Some suggestions are made for additional changes.


\end{abstract}

\subsection* {Background}

The Society for Industrial and Applied Mathematics (SIAM) publishes nearly a dozen periodicals.  \sirev\ is the ``leading'' journal which is distributed to all members.

\sirev\ underwent major revisions about a decade ago \citep {Wright1999}. Although no reason was given, the changes apparently were meant to create a publication of the highest possible profile. For example, color printing was adopted. The journal's traditional mix of papers was divided into five explicit sections:
\begin {enumerate}
\item The \sandr\ section is self-explanatory. Especially long survey papers were encouraged.
\item A \pandt\ section publishes original research papers. This section has since been renamed \expo.
\item The \sigest\ section was created to reprint important research papers from other SIAM journals.
\item The \educ\ section publishes short, instructional units.
\item \textit {Book Review\/} articles.
\end {enumerate}
When these sections were created, an older section of \pands\ was discarded. The journal also ceased printing society news in the back pages.

\subsection* {Evaluation}

\begin {figure}
\begin {center}
\includegraphics [scale=1] {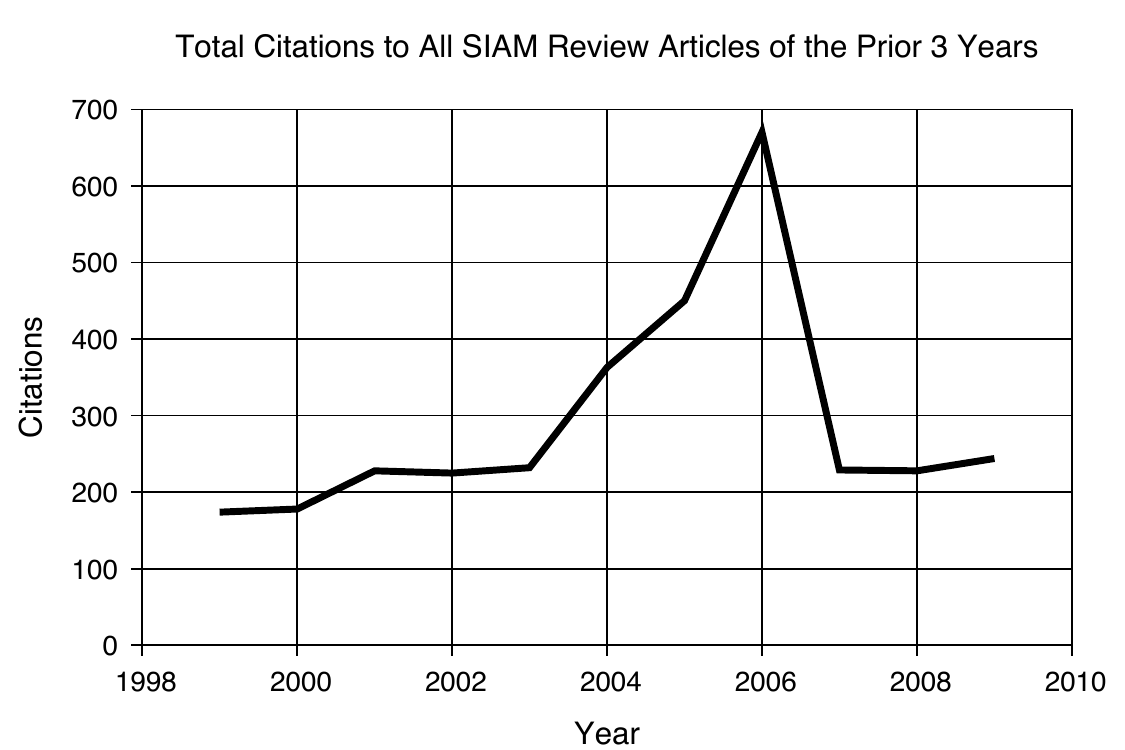}
\end {center}
\caption {Moving window of citations to \sirev\ articles of the preceding three years \citep {Scimago2007}. For example, the quantity for 2009 is the total citations then available to all articles published in 2006--2008. A single article by \citet {Newman2003} accounts for the elevated citations for the windows ending in 2004--2006.}
\label {fig:scopus}
\end {figure}

The changes to \sirev\ in 1999 did not increase the significance of the journal to the mathematical community. The three-year moving sum of citations barely rose, from 225 for 2002, to 244 for 2009 (\textbf {Figure \ref {fig:scopus}}). 

\paragraph {\sandr} The \sandr\ section has a tradition of invited submissions with sometimes spectacular results. The most highly cited article by far in any SIAM journal is an invited survey of social networks by \citet {Newman2003}. This article is off the chart compared to all others, with over 3,100 citations as of 2010. The many citations temporarily raised the five-year impact factor of \sirev, which may have wrongly suggested that the 1999 changes were successful.

\begin {figure}
\begin {center}
\hspace*{-2in} \includegraphics [scale=1] {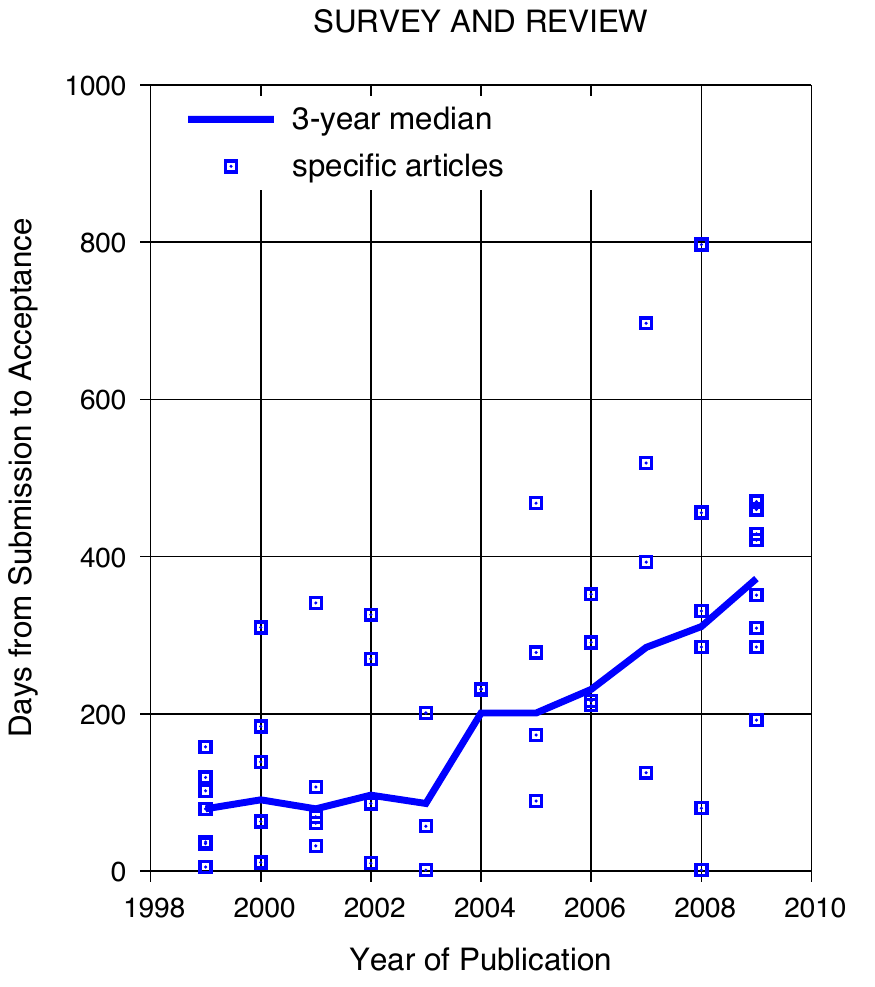}
\quad
\includegraphics [scale=1] {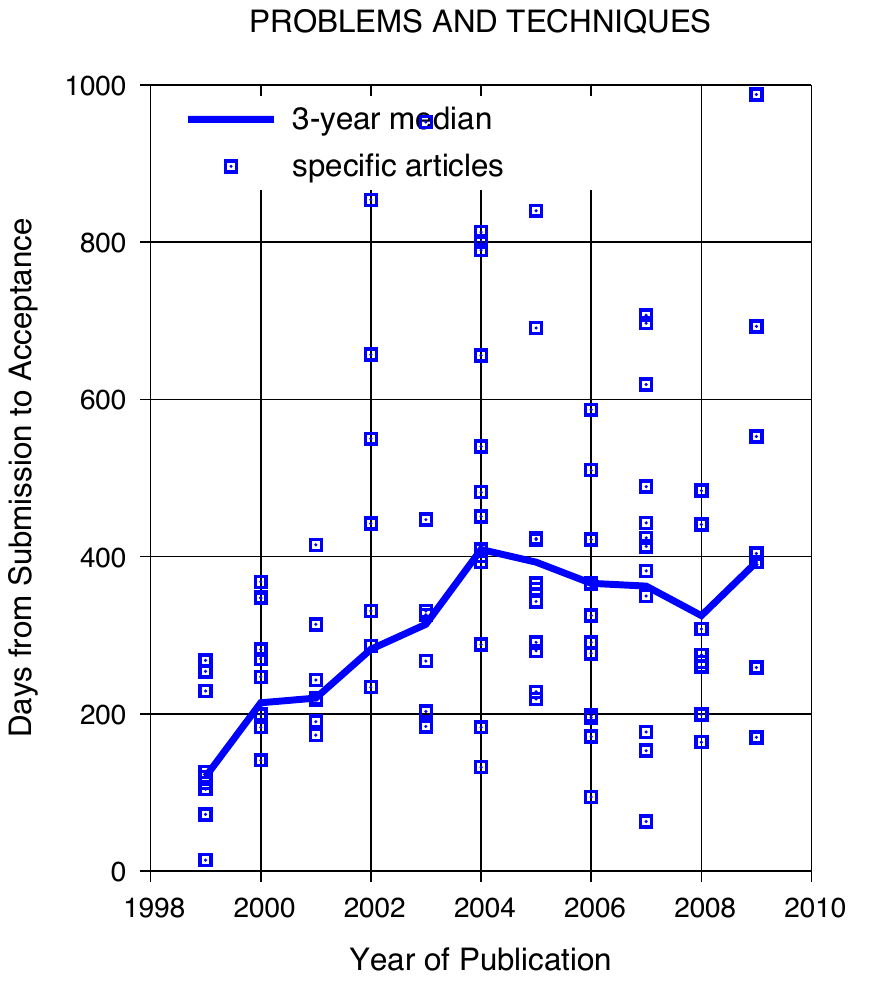}\hspace*{-2in} 
\vspace*{0.25in}

\hspace*{-2in} \hspace*{3.5in}
\quad
\includegraphics [scale=1] {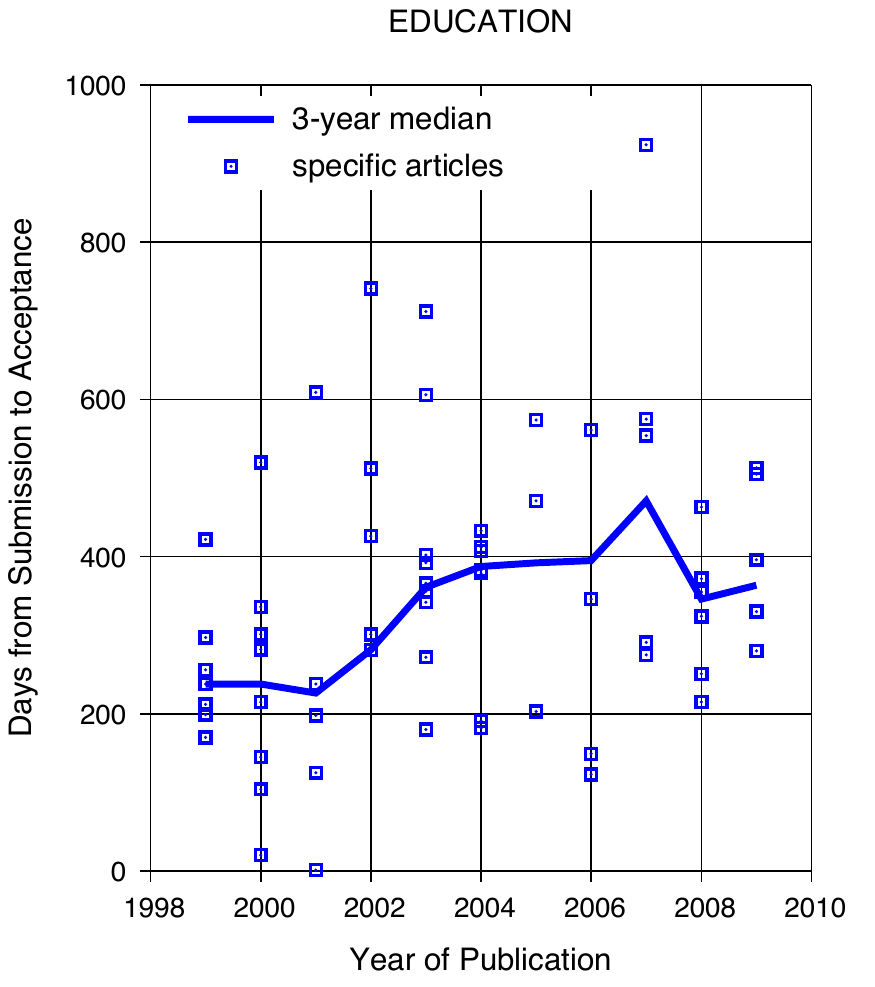}\hspace*{-2in} 
\end {center}
\caption {Duration of editorial review for articles in three sections of \sirev. Data is not shown for \sigest\ because papers in that section were refereed by other journals.}
\label {fig:days}
\end {figure}

In recent years, more \sandr\ articles appear to be contributed (as opposed to invited). This trend is suggested by the lengthening refereeing period starting in 2004 (\textbf {Figure \ref {fig:days}}, upper left).

\begin {figure}
\begin {center}
\hspace*{-2in} \includegraphics [scale=1] {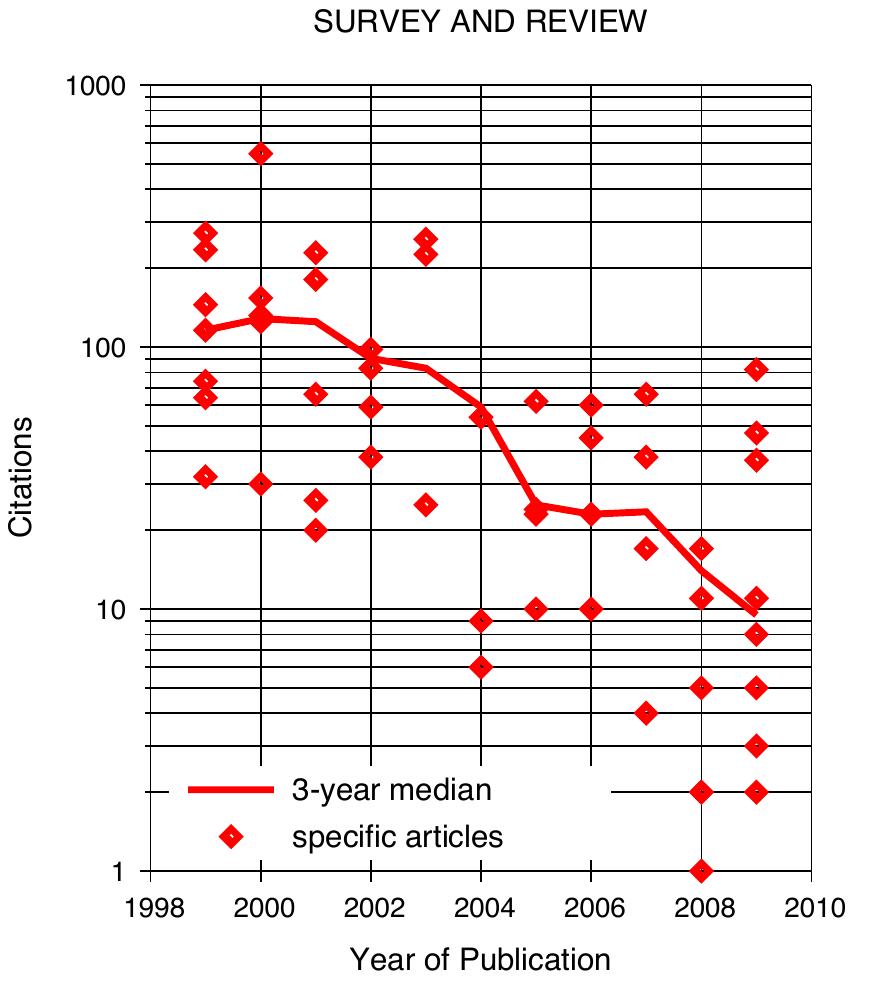}
\quad
\includegraphics [scale=1] {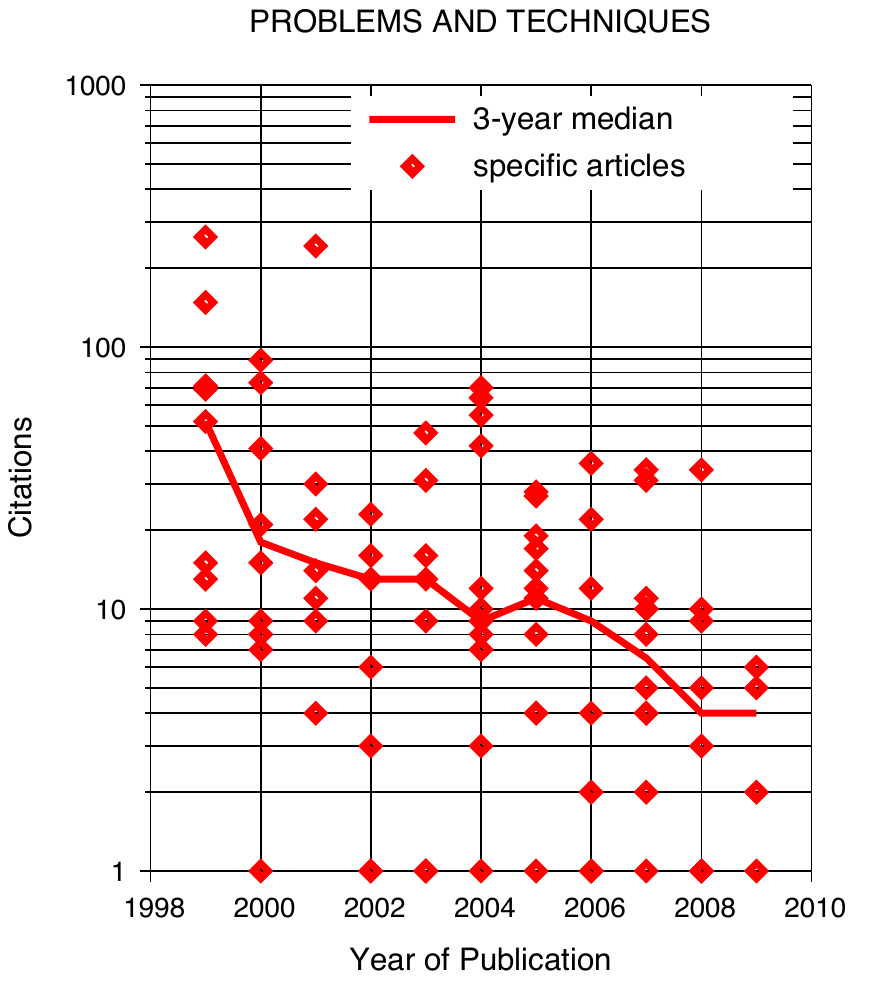}\hspace*{-2in} 
\vspace*{0.25in}

\hspace*{-2in} \includegraphics [scale=1] {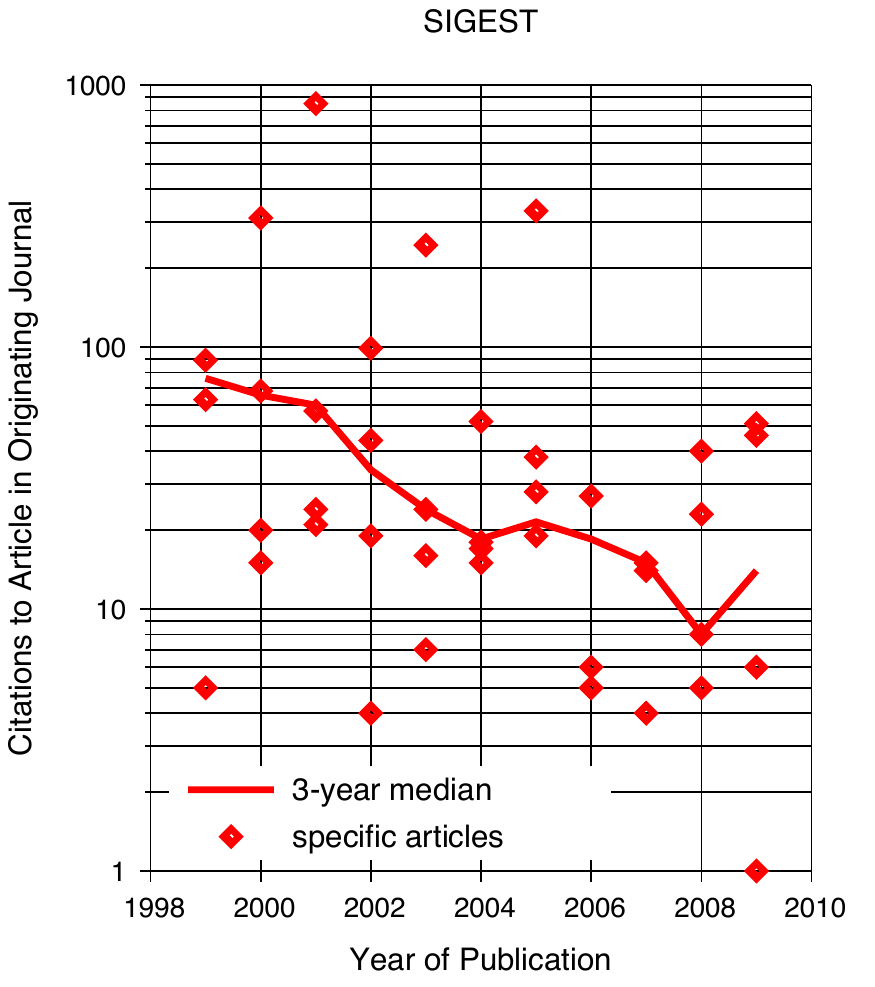}
\quad
\includegraphics [scale=1] {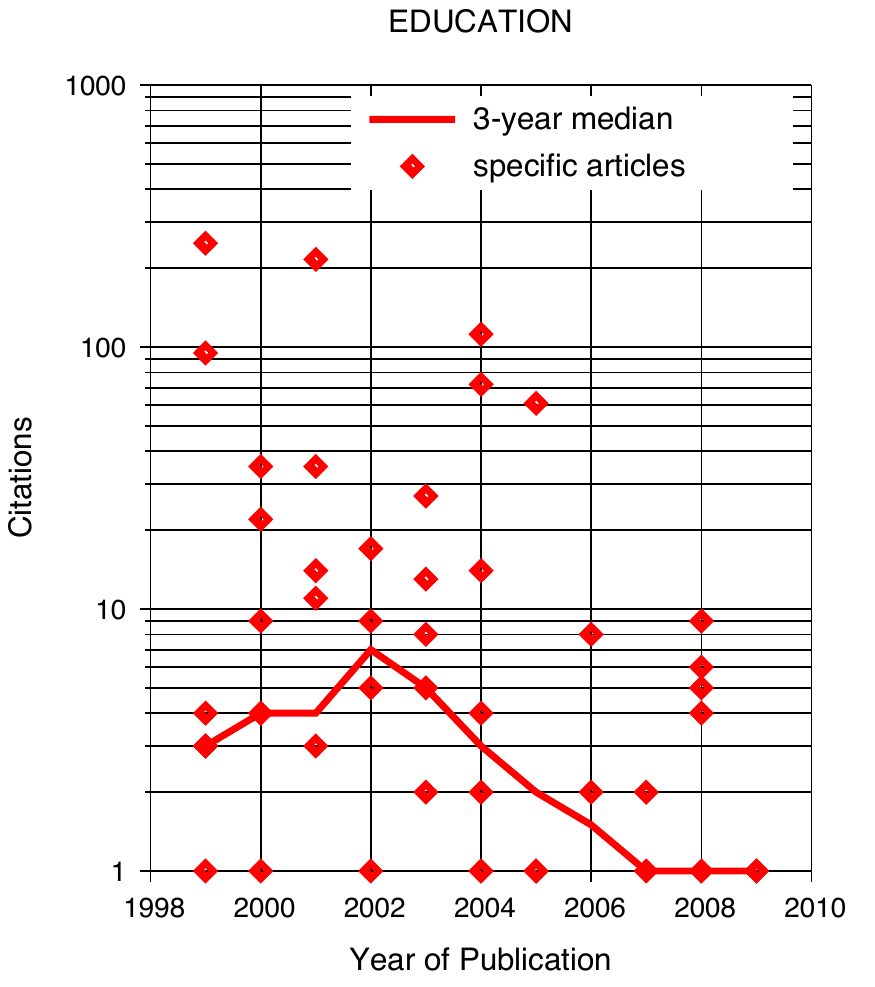}\hspace*{-2in} 
\end {center}
\caption {Quantity of citations to articles in four sections of \sirev. The most-cited \sandr\ paper of \citet {Newman2003} is off the scale.}
\label {fig:citations}
\end {figure}

Survey articles typically have more citations than research papers (\textbf {Figure \ref {fig:citations}}, upper left compared with others). Many \sandr\ articles have over 100 citations, but when the data were gathered in 2010, all papers since 2004 had fewer than 100 citations. If more \sandr\ articles are indeed contributed, then they appear to be less interesting to readers. 

\paragraph {\pandt} The \expo\ section seems not to attract important original research. The papers are consistently less frequently cited than those in \sigest\ and \sandr\ (\textbf {Figure \ref {fig:citations}}, upper right compared with left column).

\paragraph {\sigest} Some papers reprinted from other journals do receive as many citations as \sandr\ articles. It is reasonable to ask whether papers that have ten or fewer citations should be reprinted (\textbf {Figure \ref {fig:citations}}, lower left).

\paragraph {\educ}

\begin {figure}
\begin {center}
\hspace*{-2in} \includegraphics [scale=1] {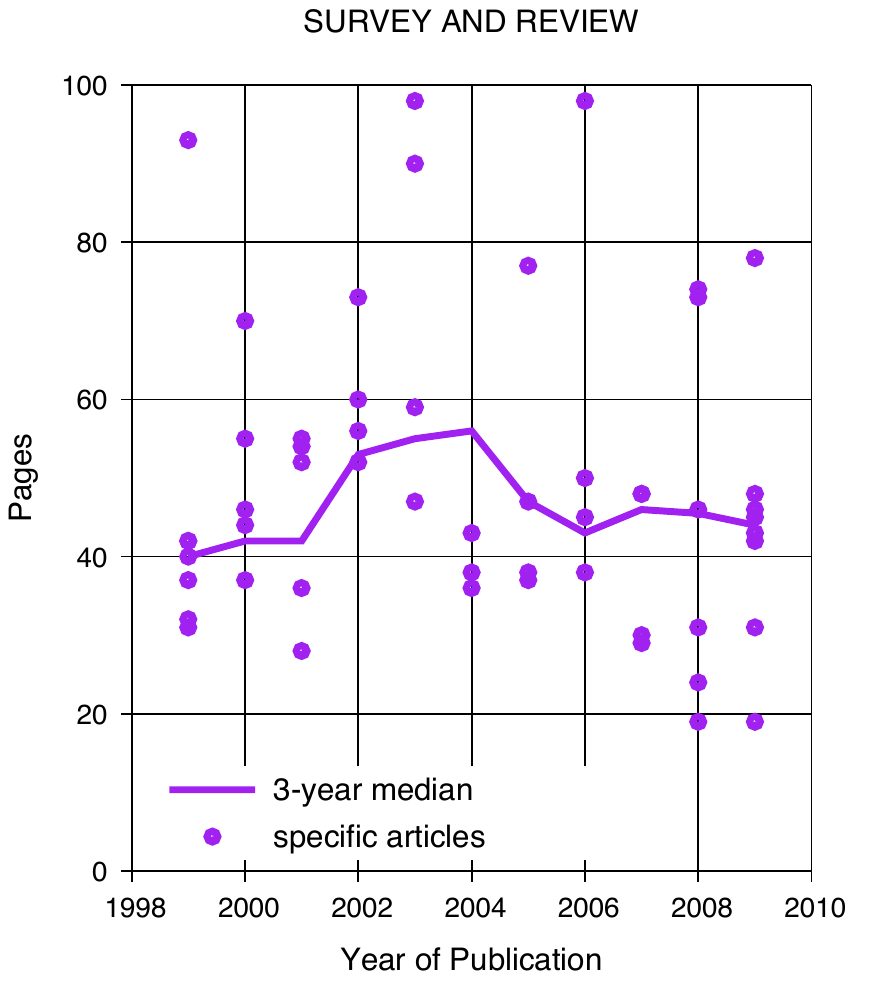}
\quad
\includegraphics [scale=1] {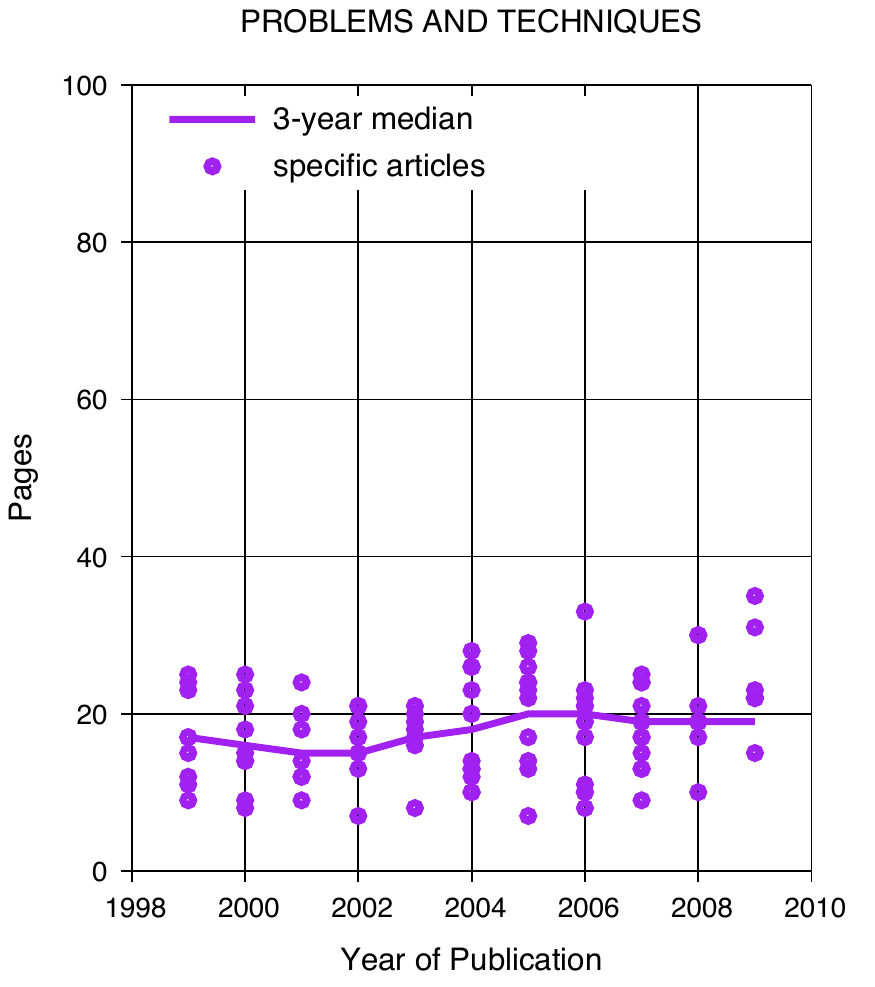}\hspace*{-2in} 
\vspace*{0.25in}

\hspace*{-2in} \includegraphics [scale=1] {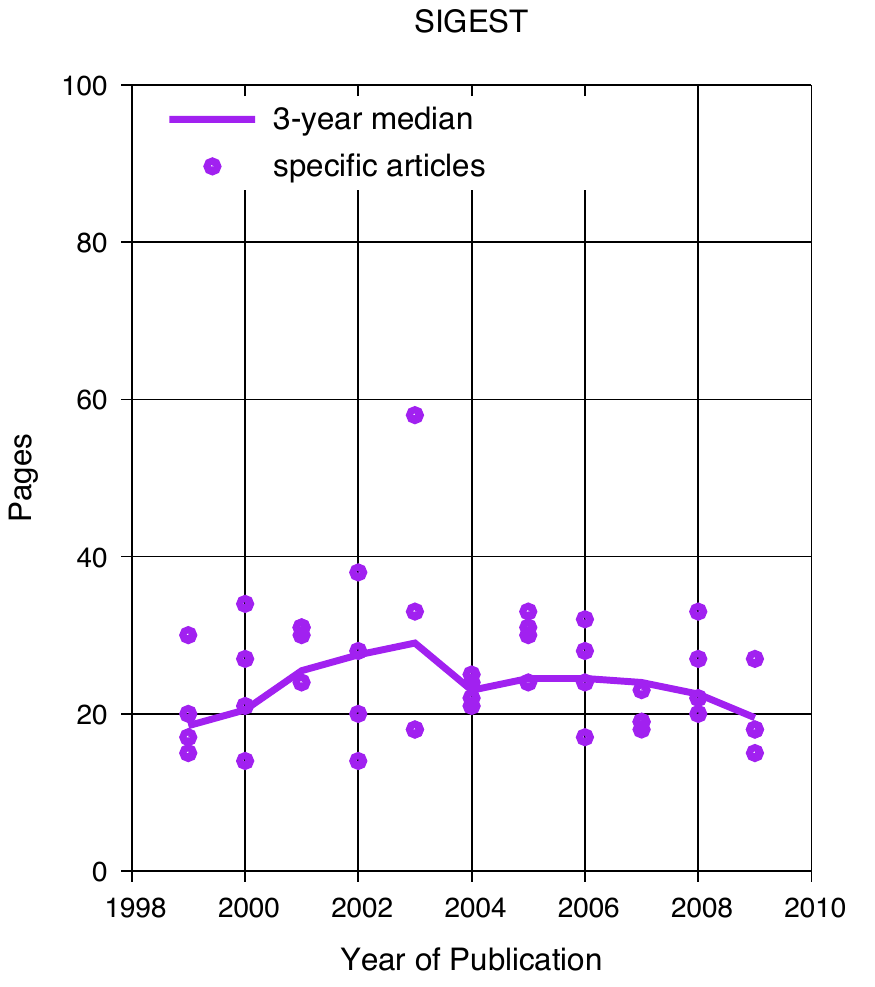}
\quad
\includegraphics [scale=1] {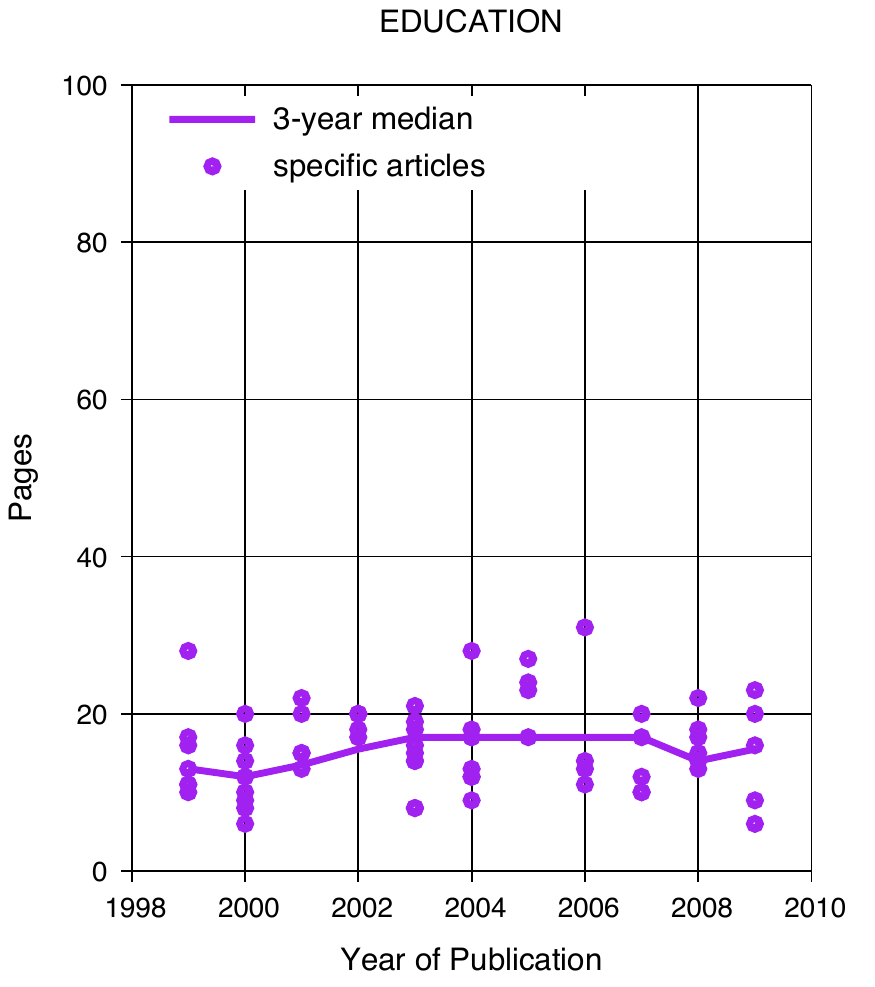}\hspace*{-2in} 
\end {center}
\caption {Length, in pages, of articles in four sections of \sirev.}
\label {fig:pages}
\end {figure}

The education papers seem to be of little interest in recent years. The median paper now receives few if any citations (\textbf {Figure \ref {fig:citations}}, lower right). One may speculate that authors are reluctant to contribute, because the refereeing period is among the longest in the journal (\textbf {Figure \ref {fig:days}}, lower right compared with upper row) even though the papers are consistently the shortest (\textbf {Figure \ref {fig:pages}}, lower right compared with others). 

\subsection* {Conclusions}

The modifications to \sirev\ a decade ago did not make the journal more valuable to readers as measured by citation analysis.  SIAM should consider further editorial changes to make the journal more attractive.  The evidence above suggest that the following changes may be beneficial:
\begin {enumerate}
\item The \sandr\ section should commission more articles that are expected to receive many citations. Independent of the editorial board, a committee of the best authors and researchers might be convened to suggest topics and authors.

\item The \expo\ section should be discontinued as there is no citation evidence that the articles are truly outstanding.  

\item \sigest\ should reprint only papers that have a great many citations. 

\item The \educ\ section should be reconsidered in light of the low citation record and the long refereeing period. 
\end {enumerate}
A purpose for modifying \sirev\ in 1999 must have been to create a stimulating publication to attract and retain society members. A publication more suited for that purpose would be a magazine rather than a research journal. With the removal of the \expo\ section, \sirev\ may well reduce the number of issues to quarterly. The savings could be used to convert \textit {SIAM News\/} to a color magazine as I have previously suggested \citep {Grcar2010g}.

\raggedright

\end{document}